\documentclass[%
 aip,
 amsmath,amssymb,
 reprint,%
]{revtex4-1}

\usepackage{graphicx}
\usepackage{dcolumn}
\usepackage{bm}

\usepackage[dvipsnames]{xcolor}
\usepackage[utf8]{inputenc}
\usepackage[T1]{fontenc}
\usepackage{mathptmx}
\usepackage{etoolbox}
\usepackage{mathbbol}
\newtheorem{theorem}{Theorem}
\makeatletter
\def\@email#1#2{%
 \endgroup
 \patchcmd{\titleblock@produce}
  {\frontmatter@RRAPformat}
  {\frontmatter@RRAPformat{\produce@RRAP{#1\href{mailto:#2}{#2}}}\frontmatter@RRAPformat}
  {}{}
}%
\makeatother
\begin{document}


\title{Chaos and Synchronization in Financial Leverages Dynamics: Modeling Systemic Risk with Coupled Unimodal Maps}



\author{Marco Ioffredi} 
\email[Corresponding author email address: ]{ioffredi@stanford.edu}
\affiliation{Stanford University, 450 Jane Stanford Way, Stanford, CA 94305, United States}

\author{Stefano Marmi}
\affiliation{Scuola Normale Superiore, Piazza dei Cavalieri, 7, 56126 Pisa PI, Italy}

\author{Matteo Tanzi}
\affiliation{King's College London, Strand, London WC2R 2LS, United Kingdom}


\date{\today}

\begin{abstract}
Systemic financial risk refers to the simultaneous failure or destabilization of multiple financial institutions, often triggered by contagion mechanisms or common exposures to shocks. In this paper, we present a dynamical model of bank leverage — the ratio of asset holdings to equity — a quantity that both reflects and drives risk dynamics. We model how banks, constrained by Value-at-Risk (VaR) regulations, adjust their leverage in response to changes in the price of a single asset, assumed to be held in fixed proportion across banks. This leverage-targeting behavior introduces a procyclical feedback loop between asset prices and leverage. In the dynamics, this can manifest as logistic-like behavior with a rich bifurcation structure across model parameters. By analyzing these coupled dynamics in both isolated and interconnected bank models, we outline a framework for understanding how systemic risk can emerge from seemingly rational micro-level behavior.
\end{abstract}


\maketitle 


\begin{quotation}
Financial systemic risk (or systemic failure) refers to the event in which several institutions fail simultaneously, possibly leading to large financial crisis. To model the evolution of the level of risk associated with a system, it turns out to be crucial to look at the evolution of banks' leverage ratios, quantities that can amplify both gains and losses.
Building on the works \cite{lillo2023unimodal}\cite{mazzarisi2019panic}\cite{corsi2016micro}, we propose a simple analytical dynamical model for the evolution of coupled leverages and explore 
how this can lead to complex behavior within the banking system. 
In particular, we investigate how the heterogeneity in bank size and strategies of single banks may impact the stability of the whole system.
These insights contribute to a deeper understanding of the mechanisms that can precipitate systemic crises and inform strategies for enhancing financial resilience.
\end{quotation}
\section{Introduction}
Several mechanisms may be responsible for the emergence of systemic risk in a system of financial institutions (or banks for simplicity), such as informational contagion leading to bank runs and liquidity shortages, direct contagion in which banks fail to pay back loans to other banks, etc.
In a broad sense, systemic risk may also be considered to include the risk associated with a system collapsing simultaneously as a result of a shock, either endogenous or exogenous (\cite{de2000systemic}). For instance, in the case in which multiple banks hold positions in the same asset, a severe loss of its value would cause all banks to fail together.

This paper will primarily focus on this broader definition of systemic risk.
As briefly explained below, the dynamics of the returns of the assets (and thus the evolution of the level of systemic risk in the ``broad'' sense) are intimately related to the evolution of the \textit{leverages} of the banks.

The leverage $\lambda$ of a market agent, such as a bank,  is defined as the ratio between the value of the assets held $A$, i.e. the economic resources the bank owns, and the value of its  equity $E$, i.e. what remains of the asset value after subtracting the liabilities $L$ to which the bank is subject to, for instance  (so that $E=A-L$). 
As the name itself suggests, a large leverage means an amplification of resources invested, which is made possible by the bank borrowing money to invest (a ``pumping'' of money from the creditors allows the bank to amplify its invested resources) and which translates into an amplification of gains, but also losses.

The relationship between leverage and risk has indeed already been suggested by many (\cite{acharya2016banks,adrian2010liquidity,adrian2014procyclical,fostel2008leverage,huizinga2012bank,nuno2017bank}), and depends on the following mechanism.
In order to maximize gain, banks want to be as leveraged as possible, as a larger leverage will in general allow for more opportunities to make profit by buying on margin (i.e. by investing borrowed money). 
On the other hand, banks have to face limitations on the maximum attainable leverage imposed by financial regulations in order to make the financial system more robust and resilient.

In particular, regulators require banks to define a risk measure, the Value at Risk (VaR) $\Lambda_{\text{VaR}}$, which is associated with a certain probability of loss $P_{\text{VaR}}$ and is implicitly defined by 
\begin{equation}
    \int_{-\infty}^{-\Lambda_\text{VaR}}f(r^p)dr^p=P_{\text{VaR}}
\end{equation} $r^p$ being the portfolio return (i.e. the
ratio between the gain/loss yielded by the investment and the
initial cost of the investment) and $f$ its probability density function.
In other words, banks are required to assure that the probability of losing an amount larger than $\Lambda_{\text{VaR}}$ is less than $P_{\text{VaR}}$.
Having defined $\Lambda_{\text{VaR}}$, it is required that the maximum loss, $A\Lambda_{\text{VaR}}$, be smaller than the equity of the bank. This translates in a maximum leverage that banks are allowed to set, i.e. $\lambda\leq \frac{1}{\Lambda_{\text{VaR}}}$.
Assume a mean $\mu$ and a variance $\sigma^2$ for $r^p$.
Then, by Chebychev inequality we can write:
\begin{equation}
\frac{\sigma^2}{(\mu+\Lambda_{\text{VaR}} )^2}\geq P(|r^p-\mu|\geq \mu+\Lambda_{\text{VaR}})\geq P(r^p\leq -\Lambda_{\text{VaR}})
\end{equation}
So that 
\begin{equation}
    P(r^p\leq -\Lambda_{\text{VaR}})\leq\frac{\sigma^2}{(\mu+\Lambda_{\text{VaR}} )^2}
\end{equation}
Requiring $\frac{\sigma^2}{(\mu+\Lambda_{\text{VaR}} )^2}\leq P_{\text{VaR}}$, i.e. $\Lambda_{\text{VaR}}\geq\frac{\sigma}{\sqrt{P_{\text{VaR}}}}-\mu$, ensures that $P(r^p\leq -\Lambda_{\text{VaR}})\leq P_{\text{VaR}}$.
We deduce the following upper bound for the leverage:
$\lambda\leq\frac{1}{\Lambda_{\text{VaR}}}\leq\frac{1}{\frac{\sigma}{\sqrt{P_{\text{VaR}}}}-\mu}$
Actually, the most stringent inequality is usually required to hold:
$\lambda\leq\frac{1}{\Lambda_{\text{VaR}}}\leq\frac{\sqrt{P_{\text{VaR}}}}{\sigma}$

Banks will thus make an estimate of $\sigma$ and set a target leverage $\lambda=\frac{1}{\alpha\sigma}$.
Having defined $\alpha=\frac{1}{\sqrt{P_{\text{VaR}}}}$.
This leverage-targeting strategy  banks may adopt can 
result in a positive feedback between leverages and asset prices (\cite{adrian2010liquidity}), possibly leading to instabilities very much like what happens in physical systems with positive feedback.
Being concrete, what happens is that if things go well (i.e. the price of the asset rises and thus $E$ increases), in order to maintain a certain leverage level, banks must increase $A$, and are thus led to borrow more (i.e. increase $L$); while if things go bad (price falls and thus $E$ decreases) banks will have to reduce $A$ (by reducing $L$). Now, the price of the asset increases (decreases) as the demand for it increases (decreases). Thus if price rise (fall), banks will increase (reduce) $A$ to reach the target leverage, i.e. the demand for the asset will rise (fall) and so will do its price, making apparent the positive feedback between asset price and leverages.

For example, consider a bank having a target leverage $\lambda=10$, which it realizes having an asset value equal to $100$ and an equity equal to $10$.
If there is a $1\%$ portfolio return, at the next time the bank will have an asset value of $101$ and an equity equal to $11$, with a resulting leverage $\frac{101}{11}\approx 9.2$, which is less than the target leverage. Thus, to increase its leverage, the bank will borrow some money and increase its asset value, determining as a consequence an increase in its price.

This positive feedback mechanism is also known as procyclicality of the Value-at-Risk constraint \cite{adrian2014procyclical} and has the effect of amplifying prices movements, which is particularly relevant during a falling period of a financial crisis (i.e. when prices of assets drop). 

On top of the above described feedback mechanism, another positive feedback effect may play a role.
Indeed, many estimates of risk are based on observation of recent price movements. However, the choice of the time window in the past to consider to perform such an estimate is far from trivial, since there is a trade off between choosing a long temporal window in order to improve statistical precision and choosing a short window in order to capture a more timely measure of risk.
This choice, determining the future trading strategies, endogenously moves the asset prices, resulting in a second feedback effect.
The models developed here will take both of these feedback effects into account.

This being said, it should be clear that by looking at the evolution of the leverages one could get insights in the level of systemic risk of the system. 

Indeed, the larger the leverages, the more amplified asset-price movements are (due to the first feedback mechanism mentioned above) and thus the riskier it is for the price to dramatically fall (maybe after a negative fluctuation), making all banks investing in it to fail simultaneously, i.e. to an increase in the level of systemic risk.

In the case in which the banks are ``similar'' in the sense that they have comparable sizes and risk estimation strategies, one can reduce to study the evolution of a ``representative'' leverage. This is what has been done for instance in a paper by Lillo et al.\cite{lillo2023unimodal}, where they obtain a slow-fast deterministic-random dynamical system for the leverage of a bank investing in a single asset.
 We take this model as a starting point, and introduce  heterogeneity in the system, allowing the banks to have different asset sizes and different risk estimation strategies.
In this case, disregarding random fluctuations, the system can be modeled with a discrete time dynamical system consisting of  unimodal maps coupled through a mean field (see equation \ref{Eq:LevEvolution} below).

In most previous works considering the impact of risk management practices on the dynamical properties of the leverages, either a single bank was considered (\cite{lillo2023analysis,aymanns2016taming}) or multiple banks and multiple assets (\cite{corsi2016micro,mazzarisi2019panic}) but with all banks being in some sense equivalent (they all solved the same optimization problem and their portfolios were determined by random choices in a pool of equivalent assets; thus, in particular, their leverages were at all instants equal).
The models introduced here try  to go towards a more realistic direction by differentiating banks on the basis of their capital size and their risk forecasting strategies.

Other works related to the dynamical modeling of systemic risk, following some different approaches, can be found in \cite{poledna2014leverage,castellacci2015modeling,choi2012financial,choi2013financial,geanakoplos2010leverage,awiszus2022market,capponi2020swing,thurner2011systemic}.
However, for example, \cite{poledna2014leverage} and \cite{thurner2011systemic} consider an agent-based model in which the agents don't determine their choices on the basis of a strategy (e.g. by solving an optimization problem) and are endowed with infinite capital, while \cite{capponi2020swing} and \cite{awiszus2022market} don't take into account the feedback mechanism due to the estimations of risk by the banks.
\\

\textbf{Organization of the paper.}
First, the case in which there is only one bank investing in a single asset (that is the case studied in \cite{lillo2023analysis,lillo2023unimodal}) is recalled.
Then, the focus will be on the case in which there are two banks investing in a single common asset. Both numerical simulations and analytical results (for some special configurations) are provided.
\section{Unimodal Evolution of the Leverage for One Bank Trading One Asset}
\begin{figure}[htp]
\centering

\includegraphics[width=\columnwidth]{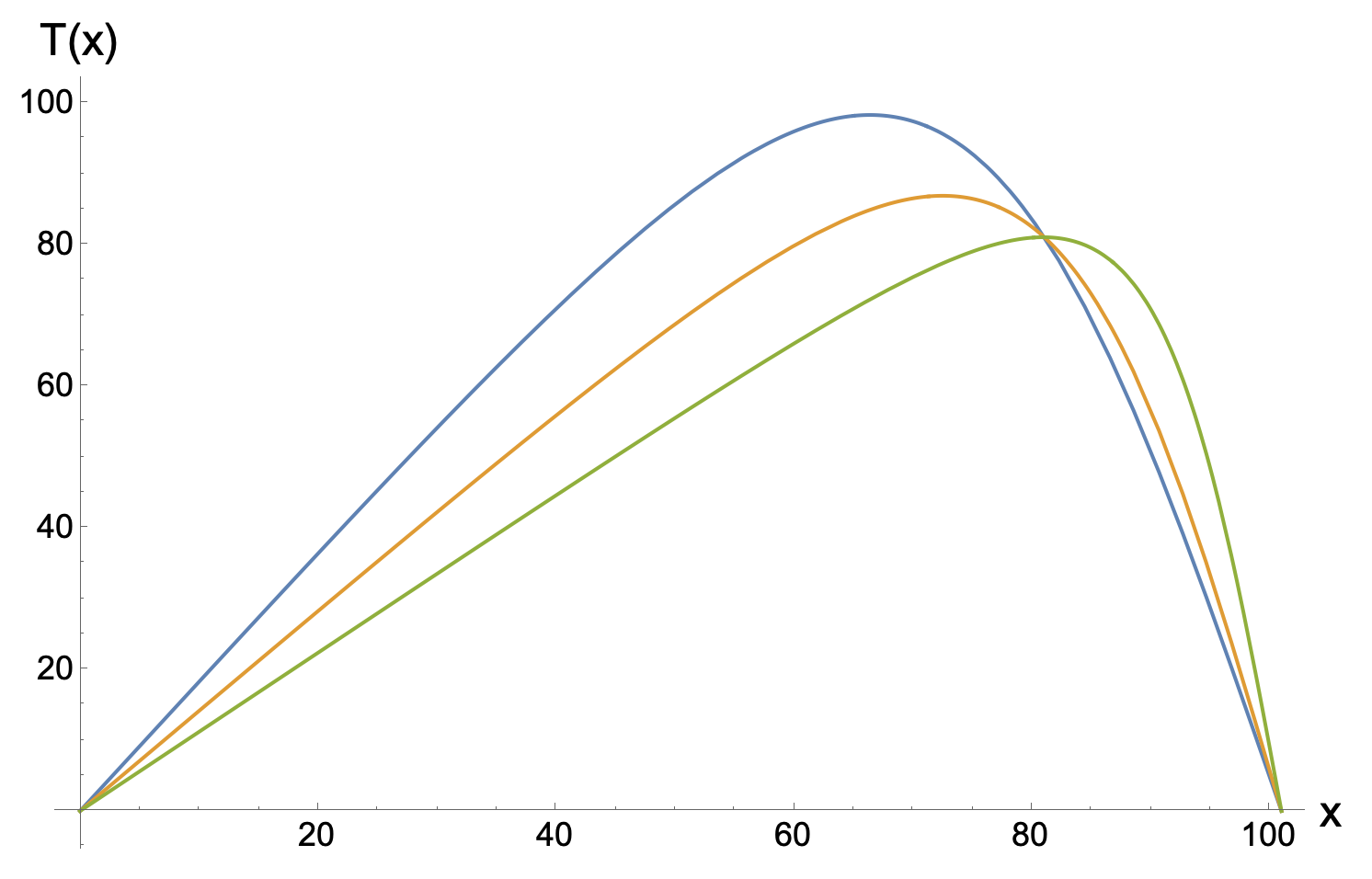}

\caption{Map $T$ for $\omega=0.3$ (Blue), $\omega=0.5$ (Orange), $\omega=0.8$ (Green). Here $\alpha=1.64$, $\Sigma_\epsilon=0.0015^2$, $\gamma=100$}
\label{fig:mapT}

\end{figure}

In the case in which there is only one bank investing in a single asset, the evolution of its leverage $\lambda_t$, with $t\in\mathbb{N}$ (in the limit in which the time scale for the financial transactions is much faster than the one used to adjust the target leverage) will turn out to be described\cite{lillo2023unimodal} by the map $T:[0,1+\gamma]\rightarrow\mathbb{R}$:
\begin{equation}
    T(x)=\left(\frac{\omega}{x^2}+\frac{(1-\omega)\alpha^2\gamma^2\Sigma_\epsilon}{(1+\gamma-x)^2}\right)^{-\frac{1}{2}}
\end{equation}
so that 
$\lambda_{t+1}=T(\lambda_{t})$.
Here $\omega\in[0,1]$ quantifies the ``memory" of the bank in forecasting the risk of the asset it is trading (a larger $\omega$ stands for a larger weight given to past observations of the asset), while $\alpha,\gamma,\Sigma_\epsilon$ are parameters whose meaning will be specified later and that we will assume to be fixed. 
The map $T$ is a unimodal map on the interval $[0,1+\gamma]$ (see Fig.\ref{fig:mapT}). 

\begin{figure}[htp]
\centering

\includegraphics[width=\columnwidth]{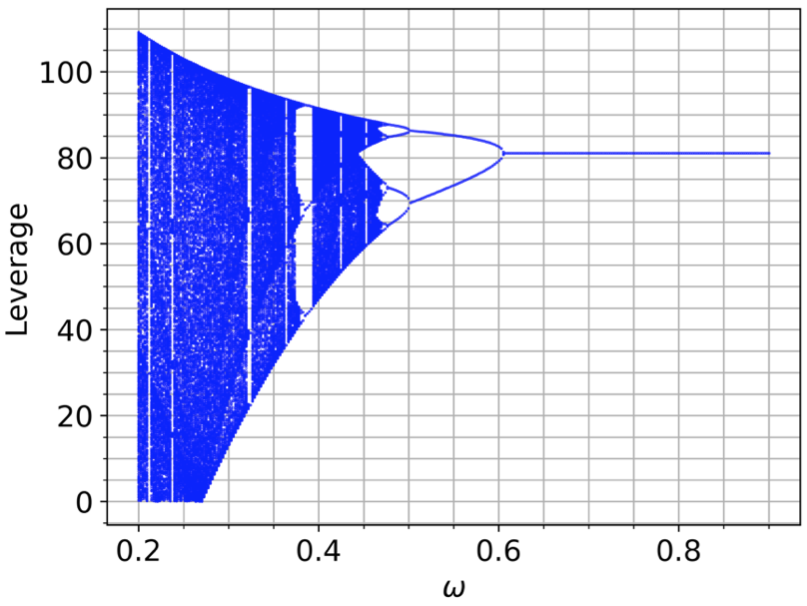}

\caption{Bifurcation diagrams for $T$ (obtained discarding the first 1000 values and considering the next 800) for $\Sigma_\epsilon=0.0015^2$, $\gamma=100$, $\alpha=1.64$.}
\label{fig:bif1}

\end{figure}

\begin{figure}[htp]
\centering

\includegraphics[width=\columnwidth]{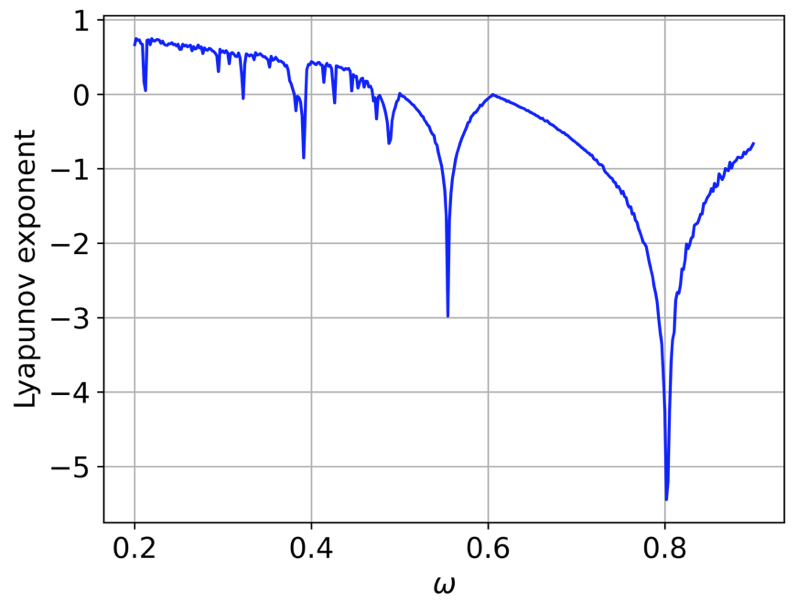}

\caption{Lyapunov exponent for the map $T$ (calculated over 100 time steps) for $\Sigma_\epsilon=0.0015^2$, $\gamma=100$, $\alpha=1.64$.}
\label{fig:lyap}

\end{figure}

By varying $\omega$, several qualitatively different behaviors are observed for the system, as shown in Fig. \ref{fig:bif1}.
Indeed, it is interesting to study how the behavior changes as $\omega$ is varied, since this is a parameter the bank can directly control.
As for the values of the other parameters employed, we follow \cite{mazzarisi2014dynamical} by setting
$\alpha = 1.64$ and $\gamma = 100$; $\Sigma_\epsilon$ has been put equal to $0.0015^2$ (thus smaller than the
one in the cited reference, as it is nevertheless to be expected since there only the
map obtained in the $\gamma\gg1$ is simulated, which is different from $T$). The choice of this value is not crucial, as the behavior of the system does not qualitatively change
for small variations of it. However, for significantly bigger values of this parameter
(e.g. $\gtrsim 0.015^2$) the behavior becomes trivial and therefore not so interesting to investigate, with the leverage reaching a fixed point for every value of $\omega\in [0, 1]$. It is also the case that, for these choices of parameters (and for $\omega$ in a subinterval of $[0, 1]$), one can find an invariant interval in $[1,1+\gamma]$ and can thus guarantee that  $\lambda_t\in[1,1+\gamma]$ for every $t\in\mathbb{N}$, as required by the interpretation of the quantity $\lambda_t$ as a financial leverage (see the supplementary material for more details on this). Lastly, these choices will be the “standard ones” in
numerical simulations in the following sections as well.
Coming back to the simulations, it is observed that for smaller values of the memory $\omega$ the system becomes more unstable, as can also be seen by looking at the Lyapunov  of the map $T$ varying $\omega$ (see Fig. \ref{fig:lyap})

Further, it can be proven that the set of $\omega$ for which $T$ is periodic is dense, and the set of $\omega$ for which $T$ is chaotic has positive Lebesgue measure (see the supplementary material for theorems useful to get these results).





\section{Coupled Unimodal Maps Model for the Leverage Evolution of Banks Trading a Common Asset}
Let's from now on focus on the  case in which there are $N$ banks investing in a single asset. This asset can be thought of as an index. 
In this case the evolution of the  leverages is given by (with $i\in\{1,...,N\}$):
\begin{equation}\label{Eq:LevEvolution}\lambda_{i,t} =\left(\frac{\omega_i}{\lambda_{i,t-1}^2}+\frac{(1-\omega_i)\alpha^2\gamma^2\Sigma_\epsilon}{\left(1+\gamma-\sum_{i=1}^2\pi_i\lambda_{i,t-1}\right)^2}\right)^{-\frac{1}{2}}
\end{equation}
Here $\pi_i\in[0,1]$ and $\sum_{i=1}^N \pi_i=1$. The $\pi$s may be given the interpretation of the ``weights'' of the banks in term of assets owned. A derivation of this model is presented in Section \ref{Sec:Model} below. But first, we are going to summarise our observations in the case $N=2$.

{\bf Remark:} \emph{We are interested in considering the behavior of the system for different values of the weights of the banks in terms of asset size $\pi_1$ and $\pi_2$, and of the ``memories'' of the banks $\omega_1$, and $\omega_2$. In other words, for our purposes the bifurcation parameters are $\omega_1,\omega_2,\pi_1$, while all others parameters are intended fixed.} 

\begin{figure}
  \centering
  
\includegraphics[width=\columnwidth]{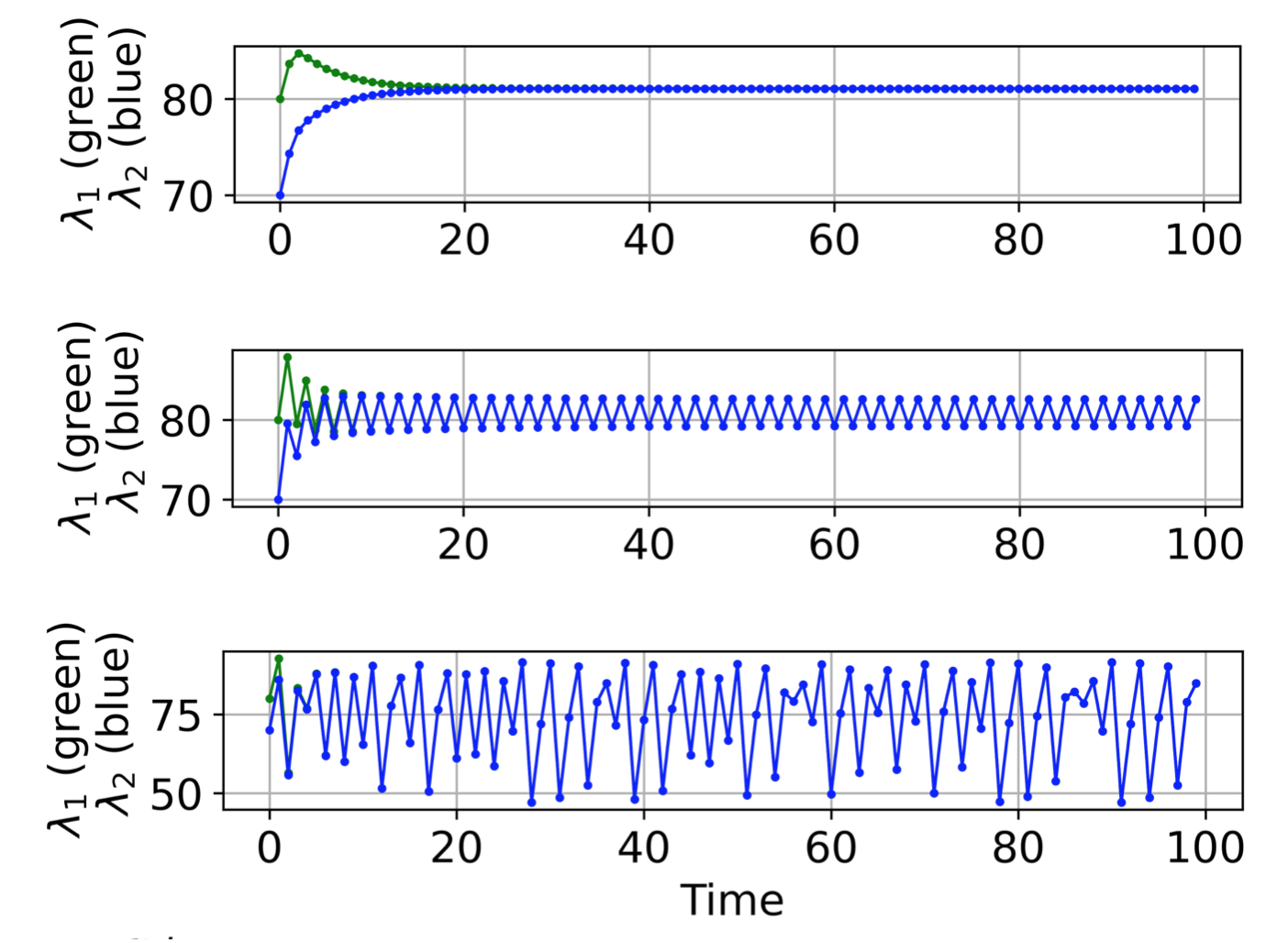}

  \caption{Leverages synchronization for $\omega_1=\omega_2=0.8, 0.6, 0.3$ (top to bottom) and $\pi_1=0.5$} 
  \label{fig:nuova}
\end{figure}

For $N=2$ and in the ``homogeneous case'' in which $\omega_1=\omega_2$, the two leverages will completely synchronize, as in Fig. \ref{fig:nuova} (see Section \ref{Sec:EqMemories}), and this happens regardless of $\pi_1$.
The common behavior of the two leverages will depend on the value of $\omega_1=\omega_2$.

{\bf Remark:}
\emph{We say that two banks are synchronised if the values of their leverages coincide. In this sense, a system of  banks is synchronized if it sits on the ``diagonal" of the phase space. In the situations where we say that the leverages synchronise, we mean that the diagonal is locally stable. An initial condition in the proximity of the diagonal will synchronize if it will approach asymptotically the diagonal or, equivalently, if the orbits of different banks asymptotically converge. Given this characterisation of synchronisation, a measurament of synchronisation could be the distance from the diagonal of the space.}

If instead the heterogeneous case is considered, in which each bank has a possibly different memory $\omega$, then the behavior can vary nontrivially depending on the values of $\omega_1$, $\omega_2$ and $\pi_1$.  

It is useful to start by considering the limiting case where bank 1 is much smaller than bank 2 to have, with good approximation,  $\pi_1=0$ (see Section \ref{Sec:Skew}). This may looked at as a forced-forcing system, where bank 2 forces bank 1.
Regardless of the initial conditions, it is observed and shown that if the forcing bank is periodic, the forced bank will also show periodic, while if the forcing bank leverage evolves chaotically, then the forced bank will behave in an irregular manner (in a sense that will be made precise later) .
Thus one may say that the behavior of the smaller bank will be determined by the larger one, whose choices thus have a crucial impact on the overall stability of the market.
This scenario generalizes to the case in which there are two groups of banks, each group with its own memory, and one of the two groups weighs significantly more than the other: in this case the memory of the group weighing more determines the nature of the evolution of the whole system.

\begin{figure}[htp]
\centering

\includegraphics[width=\columnwidth]{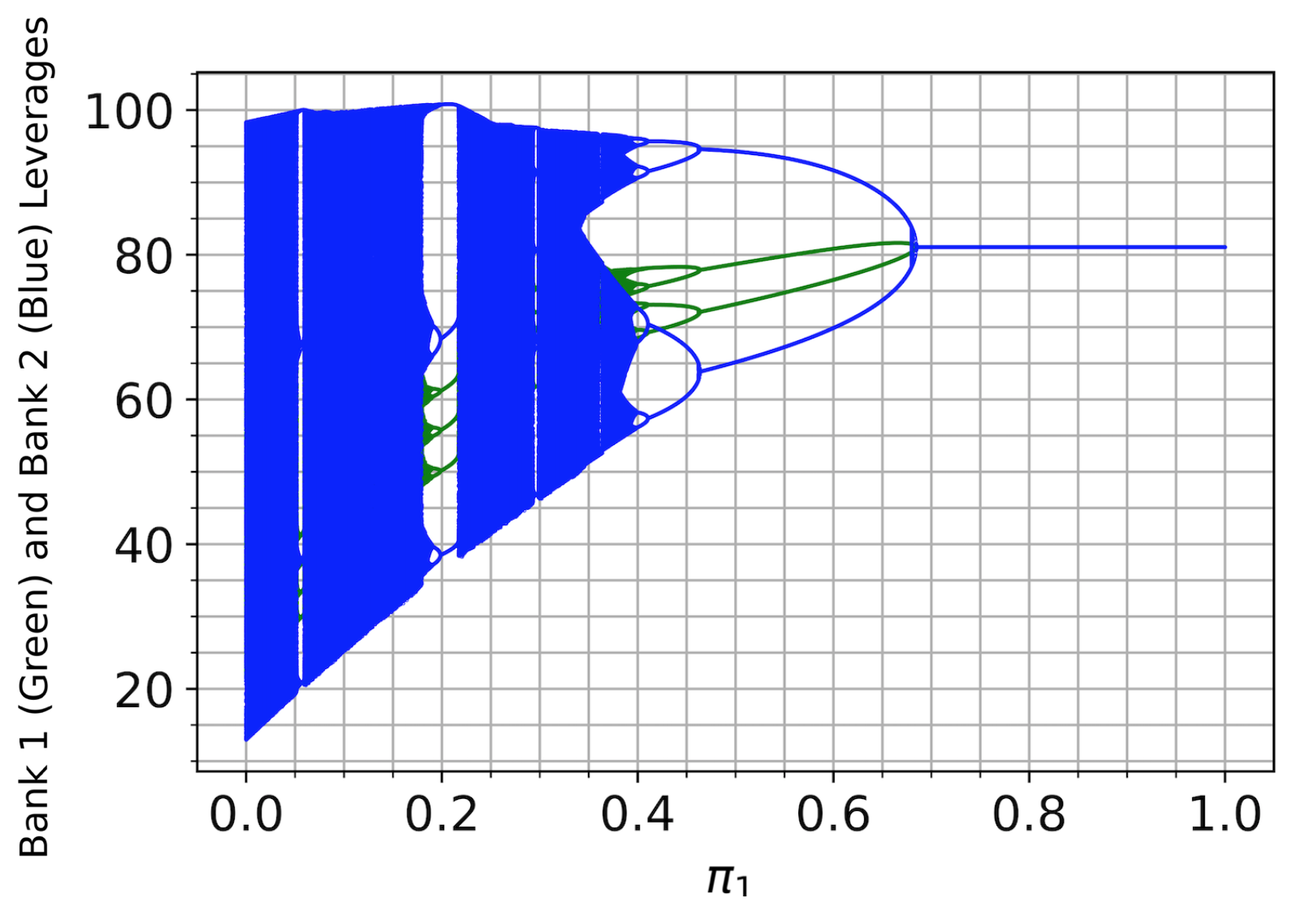}

\caption{Example of asymptotic orbits for the two banks (green for Bank 1 and blue for Bank 2) as $\pi_1$ varies and for different choices of $\omega_1, \omega_2$. Here  $\omega_1=0.8,\omega_2=0.3$. In doing the plot, the first 1000 values have been discarded and the next 500 plotted.}
\label{fig:bif2_1}

\end{figure}

Considering now the case in which $\pi_1\neq0,1$, interesting bifurcation diagrams for the leverages of the two banks are observed as relevant parameters are changed.
In particular, by keeping $\omega_1$ and $\omega_2$ fixed and changing the weight $\pi_1$ one observes a bifurcation diagram which  ``connects'' the behaviors one would have in the forced-forcing systems corresponding to $\pi_1=0$ (where only $\omega_2$ has a role in determining the nature of the behavior) and $\pi_1=1$ (where now $\omega_1$ is the relevant parameter).
Therefore, for example, there exist values for $\omega_1$, $\omega_2$ for which there is a threshold value for $\pi_1$ below which both banks switch from a regular behavior to an aperiodic one (see Fig. \ref{fig:bif2_1})

\begin{figure}[htp]
\centering

\includegraphics[width=\columnwidth]{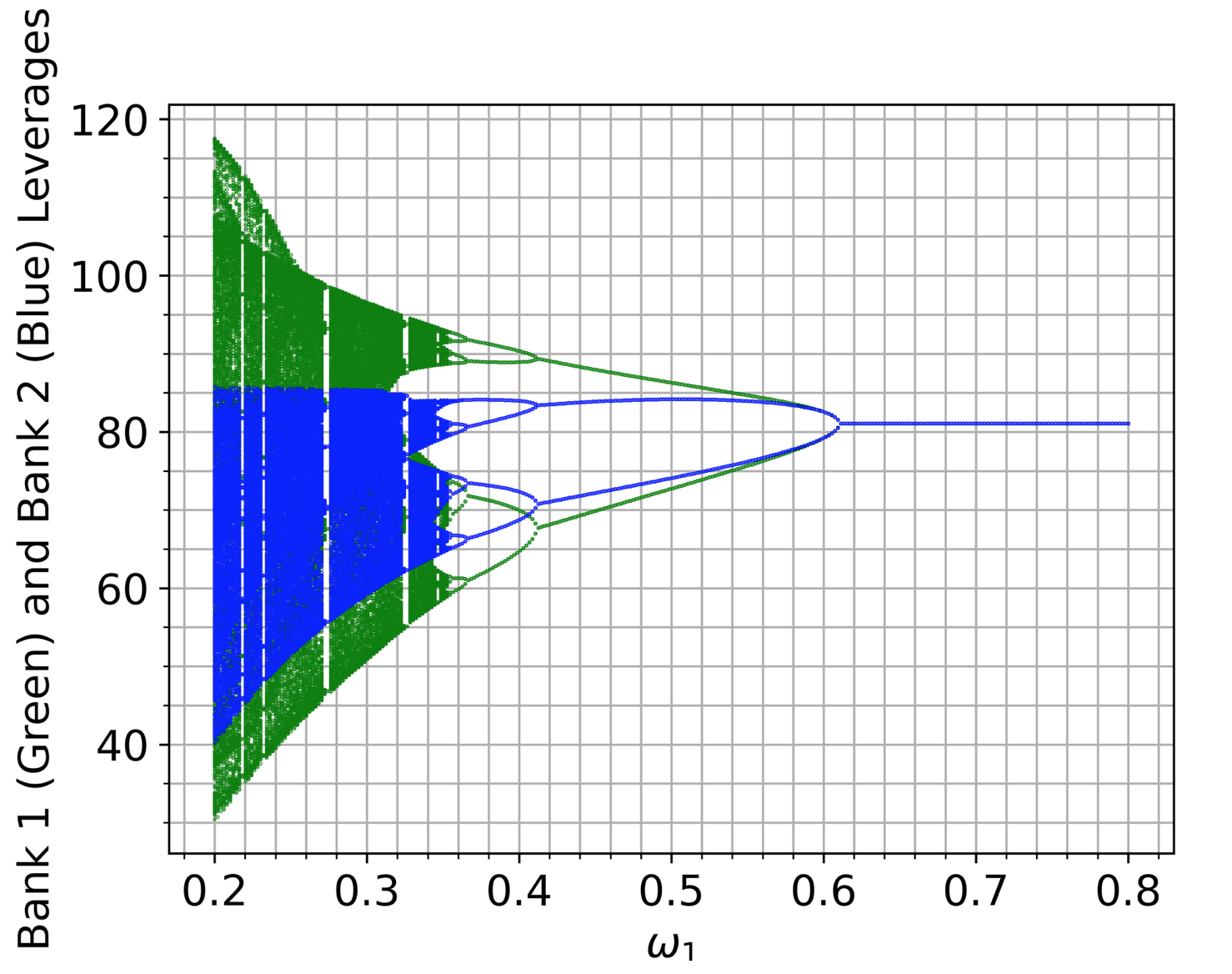}

\caption{Example of asymptotic orbits for the two banks (green for Bank 1 and blue for Bank 2) as  $\omega_1$ varies for $\pi_1=0.5$ and $\omega_2=0.6$. In doing the plot, the first 1000 values have been discarded and the next 800 plotted.}
\label{fig:bif2_2}
\end{figure}

\begin{figure}[htp]
\centering

\includegraphics[width=\columnwidth]{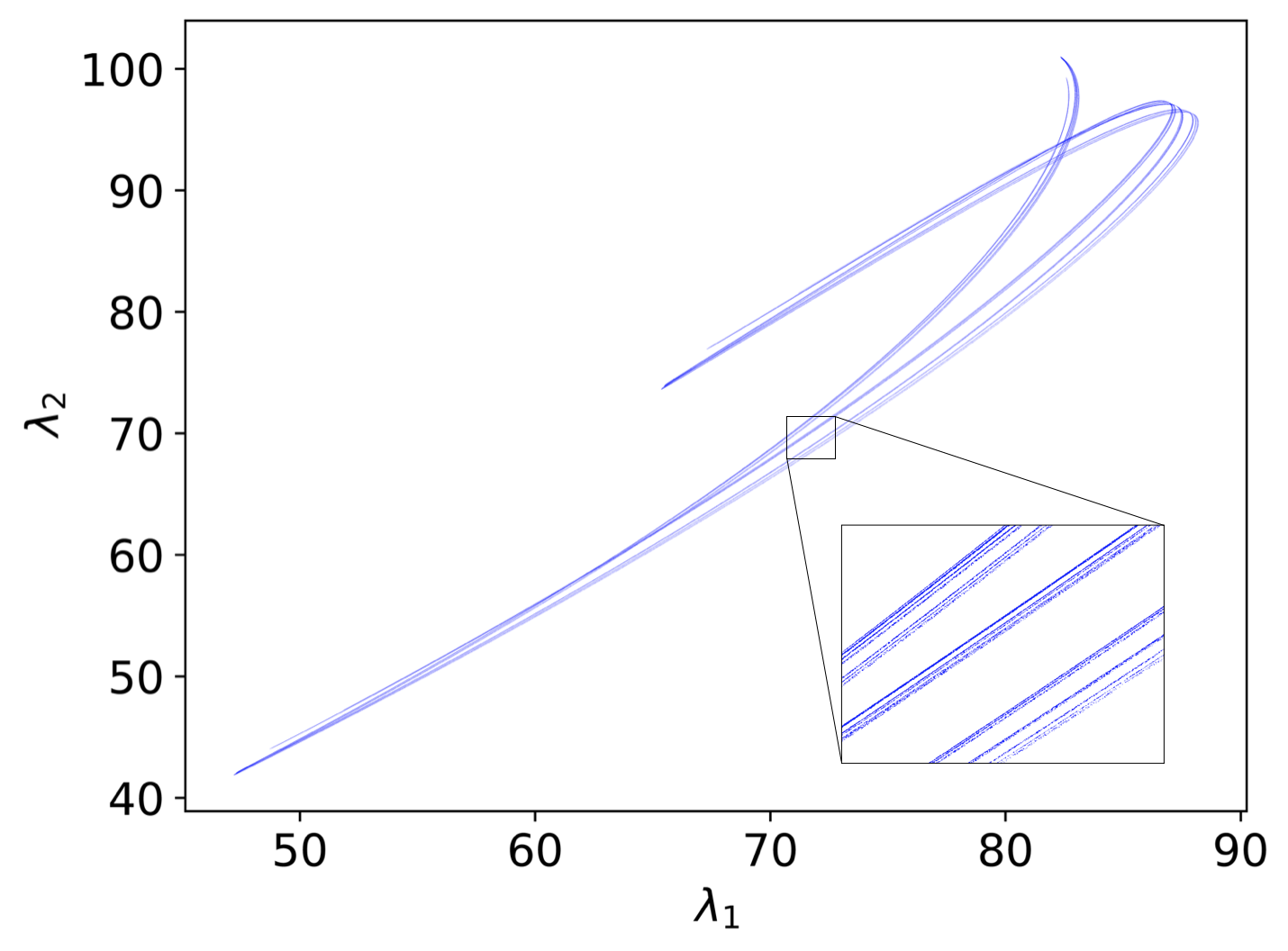}

\caption{Hénon-like attractor in the $\lambda_1$ vs $\lambda_2$ plane. Here $\pi_1=0.5,\omega_1=0.5,\omega_2=0.3$
}
\label{fig:henon}

\end{figure}

\begin{figure}[htp]
\centering

\includegraphics[width=\columnwidth]{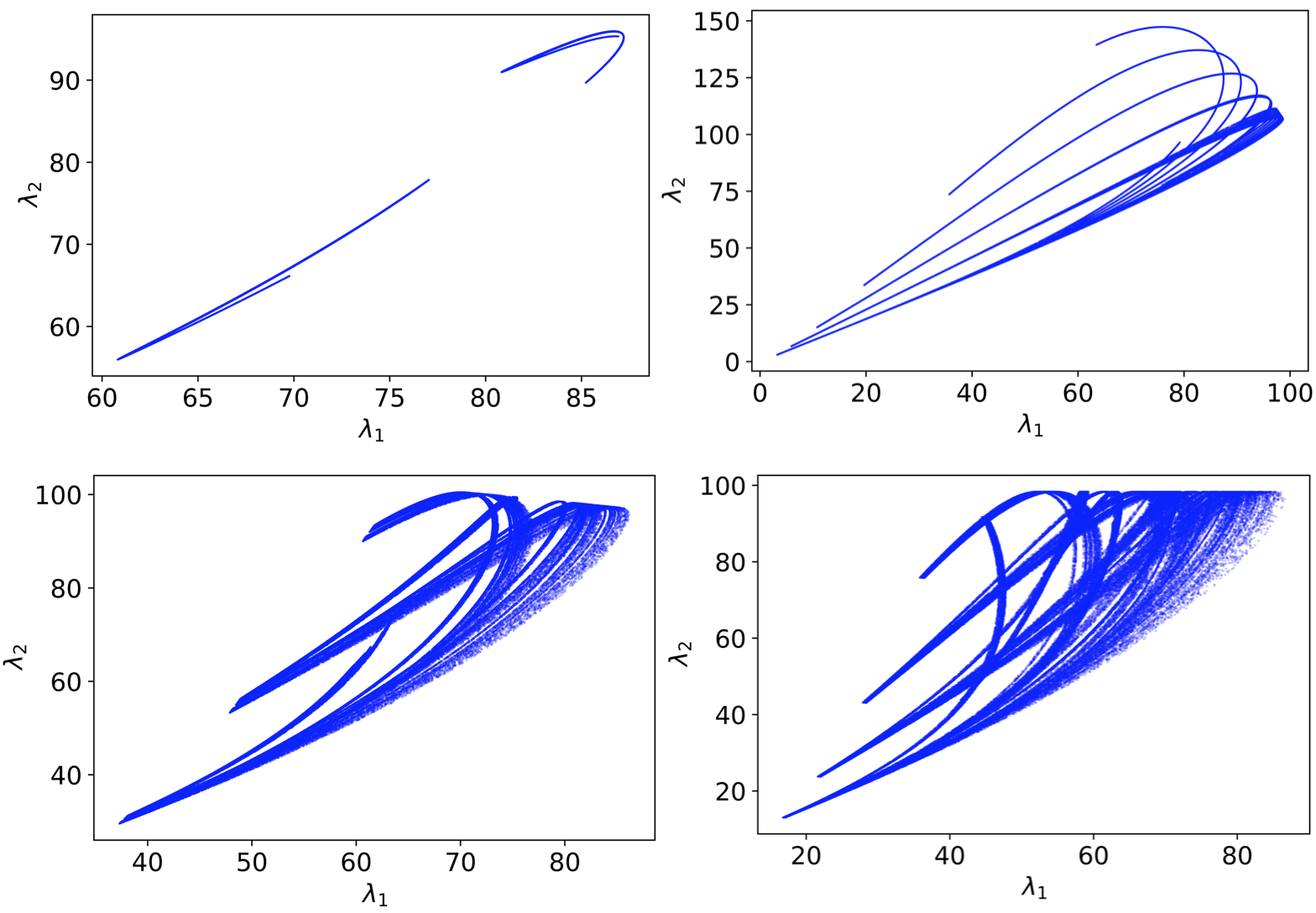}

\caption{Attractor in the $\lambda_1$ vs $\lambda_2$ plane. From top to bottom, left to right, $\pi_1=0.8,\omega_1=0.5,\omega_2=0.3$; $\pi_1=0.8,\omega_1=0.3,\omega_2=0.2$; $\pi_1=0.2,\omega_1=0.6,\omega_2=0.3$; $\pi_1=0.001,\omega_1=0.6,\omega_2=0.3$
}
\label{fig:attractors}

\end{figure}

One could also look at what happens if $\pi_1$ is fixed and $\omega_1, \omega_2$ are varied. This may be of particular interest as the memories can be more easily and directly controlled by the bank than the weights.
In this case, having fixed $\pi_1$, one observes for example that for every $\omega_2$ above a certain threshold, there is a critical $\omega_1$ value above which the system is stable (i.e. both leverages reach a common fixed point) and below which one observes period doubling bifurcations leading to chaos (see Fig. \ref{fig:bif2_2}). This critical value depends on $\pi_1$ in such a manner that as $\pi_1$ gets closer to one it becomes less and less dependent on $\omega_2$. Intuitively, for large $\pi_1$ the stability of the system depends almost exclusively on $\omega_1$. More plots from numerical simulations can be found in the supplementary material.
Next, it is also interesting to note that whenever the behavior of the two banks is aperiodic, there is not a functional asymptotic dependence between the two leverages.
Instead, for some choice of parameters, one observes in the $\lambda_1,\lambda_2$ plane an Hénon-like attractor (see Fig. \ref{fig:henon} and the supplementary material). Indeed this may not be so surprising as we are dealing with a system constituted of perturbed unimodal maps which, in certain limits, can indeed be looked at as a linear transformation of an Hénon system.
It is indeed found that the attractor has a nontrivial box counting dimension of $\approx 1.2$.
For other choices of parameters, one observes other kind of attractors, as shown for instance in Fig. \ref{fig:attractors}
\subsection{The Model}\label{Sec:Model}
In the case of a single bank, referring to the model introduced in \cite{corsi2016micro} and further developed in \cite{lillo2023unimodal,mazzarisi2019panic},
the evolution of its leverage will be the result of the following main principles: 1) to increase profit,
a bank wants to maximize its leverage up to the limit imposed by regulations; 2) this limit depends on the risk associated to the asset that is estimated by the bank in a way that can vary from bank to bank; 3) the time evolution of the prices of the asset is characterized by an autoregressive process with a time scale faster than the one in which the banks update their leverage; 4) the relative size of the total assets detained by the banks (i.e. their weights $\pi_i$) are approximately constant. 

{The last assumption simplifies the model making it tractable analytically and leading to the theorems below. Although not universally verified, we expect it to hold in the situations we study analytically: when banks have synchronized leverages due to the same risk estimation strategies, and when one bank is several orders of magnitude larger than the others (i.e. $\pi_1\approx 1$). We also expect it to hold more generally in the short term, or until large shocks perturb the system (which goes beyond the scope of our model).}

This results in a slow-fast deterministic-random dynamical system with heteroschedastic noise. In the limit in which the fast time scale is much faster than the slow one, this reduces to a deterministic dynamical system given in \eqref{Eq:LevEvolution}. 
We now show how to derive this expression by generalising the derivation given in \cite{lillo2023analysis} to the case of multiple banks. 

At each time $t\in\mathbb{N}$, bank $i$ has an equity $E_{i,t}$ and an amount of the asset, $A_{i,t}$. The ratio $\lambda_{i,t}=\frac{A_{i,t}}{E_{i,t}}\geq 1$ is called the leverage of the bank.
Each bank will try to maximize its leverage (in order to increase gains), but it will have to face VaR type constraints imposed by financial institutions.
The dynamics of the leverages is determined by two main interactions between the banks and the asset prices. First, VaR constraint determine the leverage as a function of the behavior of asset prices.
Second, a given target leverage will determine a sequence of trading events which will impact the price of the asset.
Thus one sees that it is possible to ``close the circle'' and obtain a law for the time evolution of the leverages.

If $\sigma_{e,t}$ is an estimate of portfolio variance made at time $t$, VaR constraints require that the leverage of bank $i$ should be such that $\alpha\sigma_{e,t}A_{i,t}\leq E_{i,t}$, where $\alpha$ depends on the distribution of the returns of the asset and on the strictness of the constraints (for instance, it is $1.64$ for Gaussian returns and for a VaR probability of $5\%$).
Thus one may assume that
\begin{equation}\label{eq:fond}
\lambda_{i,t}=\frac{1}{\alpha\sigma_{e,t}}
\end{equation}
Now consider how fixing a target leverage impacts the movement of the asset's price.
Having set a value for the target leverage, the bank will trade the asset between time $t$ and $t+1$ in order to keep its leverage equal to this target value.
This trading process takes places at a faster pace than the one in which the target leverage is updated (i.e. at every integer time). Let's thus introduce a quantity $\mathbb{n}\in\mathbb{N}$ so that trading operations occur $\mathbb{n}$ times every integer time step.
Next, call $r_{t+\frac{k}{\mathbb{n}}}$, with $k=1,2,\dots,\mathbb{n}$ the returns on investments (as already defined in the introduction). The dynamics of these returns can be thought of as made of two components:
\begin{equation}
    r_{t+\frac{k}{\mathbb{n}}}=\epsilon_{t+\frac{k}{\mathbb{n}}}+e_{t+\frac{k-1}{\mathbb{n}}}
\end{equation}
the first term in the RHS is an exogenous component (depending on external, non modeled events) given by a white noise term with variance $\sigma^2_\epsilon$, while the second term in the RHS is instead an endogenous component, which is the one depending on the trading actions of the bank. More in detail, the endogenous component $e_{t+h/\mathbb{n}}$ depends on the  demand for the asset arising from portfolio rebalancing of the bank, as will now be explained. 
Given a fractional time $s=t+\frac{h}{\mathbb{n}}, h=1,2,\dots,\mathbb{n}$, let's call the desired asset size for a generic bank $i$ $A^\star_{i,s}\doteq \lambda_{i,t}E_{i,s}$. Therefore at $s$ bank $i$ rebalances its portfolio by trading (i.e. buying or selling, depending on the sign of $A^\star_{i,s}-A_{i,s}$) the quantity

    \begin{align}
    A^\star_{i,s}-A_{i,s}&=\lambda_{i,t}E_{i,s}-A^\star_{i,s-\frac{1}{\mathbb{n}}}(1+r_{s})\\ &=\lambda_{i,t}(E_{i,s-\frac{1}{\mathbb{n}}}+r_{s} A^\star_{i,s-\frac{1}{\mathbb{n}}})-A^\star_{i,s-\frac{1}{\mathbb{n}}}(1+r_{s})\\ &=(\lambda_{i,t}-1)r_{s} A^\star_{i,s-\frac{1}{\mathbb{n}}}
    \end{align}

Where it is assumed that from $t$ to $t+1$ the liabilities of the banks are constant, so that the change in their total assets derives only from that of the equity.
The total demand for the asset will then be given by \begin{equation}
    D_{s}=\sum_{i=1}^N (A^\star_{i,s}-A_{i,s})
\end{equation}

Assuming a standard linear price impact function (i.e. assuming that the change in the price of the asset is directly proportional to the demand for it), the endogenous component $e_{s}$ for the return of the asset is given by \begin{equation}
    e_{s}=\frac{D_{s}}{\gamma C_{s}}
\end{equation} with $\gamma$ quantifying the liquidity of the asset and where \begin{equation}
    C_{s}=\sum_{i=1}^N A^\star_{i,s-\frac{1}{\mathbb{n}}}
\end{equation} is a proxy for market capitalization of the investment. 
What we are saying here is that the price change is larger when the ratio of the demand for it to the total amount of money invested in the asset is large and when the liquidity is low (indeed, by definition the liquidity measures the impact of selling/buying an asset on its price).

Thus:
\begin{equation}
    r_s=\frac{\sum_{i=1}^N (\lambda_{i,t}-1)A^\star_{i,s-\frac{2}{\mathbb{n}}}}{\gamma \sum_{i=1}^N A^\star_{i,s-\frac{2}{\mathbb{n}}}}r_{s-\frac{1}{\mathbb{n}}}+\epsilon_s
\end{equation}

Let's put, for every $i\in\{1,2,\dots,N\}$, $\pi_{i,s}=\frac{ A_{i,s}^\star}{\sum_{a=1}^N  A_{a,s}^\star}\in[0,1]$, so that $\sum_{i=1}^N\pi_{i,s}=1$.

Now notice that (recall $A^\star_{i,s}\doteq \lambda_{i,t}E_{i,s}$ and $E_{i,s}=\prod_{q=0}^s(1+r_q)E_{i,s}$) if the leverages were the same for all banks, then $\pi_{i,s}$ would not depend on time; this is a consequence of the fact that all banks are investing in the same asset. More precisely, the weights in this case would be constant both within a unit time interval (e.g. between $t-1$ and $t$) and across different times $t$. Thus, in the homogeneous case we can consider the weights as independent of $s$. Moreover, numerical simulations show that the error made by considering the weights as independent of $s$ even in the heterogeneous case is small enough to make it reasonable to start by analyzing the simplest case in which the weights are constant in time, leaving the investigation of the general case in which the weights change over time to future works.
All the more reason for the fact that for the cases studied analytically more in depth (i.e. the homogeneous case and the forced-forcing one) this approximation is exact.
Let's therefore put $\pi_{i,s}=\pi_i$ for every $s\in\frac{\mathbb{N}}{\mathbb{n}}$ and look at $\pi_i$ as a measure of the ``size'' of bank $i$.
It is then possible to write:
\begin{equation}
    r_s=\frac{\sum_{i=1}^N(\lambda_{i,t}\pi_i-1)}{\gamma}r_{s-\frac{1}{\mathbb{n}}}+\epsilon_s
\end{equation}

So that the dynamics of the returns of the asset for times in $(t,t+1]$ may be modeled again as a AR(1) process with auto regressive parameter \begin{equation}
\phi_t=\frac{\sum_{i=1}^N(\lambda_{i,t}\pi_i-1)}{\gamma}
\end{equation} which now depends on a convex combination of the leverages of the $N$ banks $\lambda_1, \lambda_2, \dots, \lambda_N$ in a mean field fashion.
At time $t+1$, banks need to update their target leverage so that $\lambda_{i,t+1}=\frac{1}{\alpha\sigma_{e,t+1}}$.
Thus, it now remains to provide an expression for the forecast of the variance of the aggregate return of the asset $\sigma_{e,t+1}$.
A reasonable hypothesis is to assume that the forecast made by bank $i$ at time $t+1$ is a weighted average ($\omega_i\in[0,1]$ being the weight) of the previously made forecast of the same kind and a statistical estimate of aggregate returns made by observing the returns between $t$ and $t+1$. I.e., for bank $i$:
\begin{equation}\label{eq:sigma}\sigma_{e,t+1}^2=\omega_i\sigma_{e,t}^2+(1-\omega_i)\hat{\sigma}^2_{e,t+1}
\end{equation}
and \begin{align}
    \hat{\sigma}^2_{e,t+1}&=\widehat{\text{Var}}\left[\sum_{k=1}^\mathbb{n}r_{t+\frac{k}{\mathbb{n}}}\right]
\end{align}
which is the aggregated variance of the AR(1) taking place between $t$ and $t+1$ as a function of the estimates of the parameters $\phi_{t}$ and $\sigma_\epsilon^2$, namely $\hat{\phi}_{t}$ and $\hat{\sigma}^2_\epsilon$.
In the $\mathbb{n}\rightarrow\infty $ limit, one has $\hat{\sigma}_{e,t+1}^2=\frac{\mathbb{n}\hat{\sigma}^2_\epsilon}{(1-\hat{\phi}_{t})^2}$ (see Appendix A of \cite{mazzarisi2019panic} for the derivation of this result), having introduced the Maximum Likelihood estimates for $\sigma_\varepsilon$ and $\phi_{t}$.
These are given by (see again \cite{mazzarisi2019panic}):
\begin{equation}
    \hat{\phi}_t=\frac{\sum_{k=1}^\mathbb{n}r_{t+k/\mathbb{n}}r_{t+(k-1)/\mathbb{n}}}{\sum_{k=1}^\mathbb{n}r^2_{t+(k-1)/\mathbb{n}}}
\end{equation}
and
\begin{equation}
    \hat{\sigma}_\epsilon^2=\frac{\sum_{k=1}^\mathbb{n}(r_{t+k/\mathbb{n}}-\hat{\phi}_{t}r_{t+(k-1)/\mathbb{n}})^2}{\mathbb{n}}
\end{equation}
Moreover, in the $\mathbb{n}\rightarrow\infty $ limit one expects the limit $\lim_{\mathbb{n}\rightarrow\infty}\mathbb{n}\sigma_\epsilon^2\doteq \Sigma_\epsilon$ to exists. Indeed, as already noticed one may consider the AR(1) as the discretization (with discretization step $\frac{1}{\mathbb{n}}$) of an Orstein-Uhlenbeck process, and so by a scaling argument the existence of the above limit may be deduced.
Combining equations \ref{eq:fond} and \ref{eq:sigma}  one thus has, in the $\mathbb{n}\rightarrow\infty $ limit,
(for $i\in\{1,2,\dots,N\}$)
\begin{equation}\label{eq:2}
\lambda_{i,t+1}=\left(\frac{\omega_i}{\lambda_{i,t}^2}+\frac{(1-\omega_i)\alpha^2\Sigma_\epsilon}{(1-\hat{\phi}_{t})^2}\right)^{-\frac{1}{2}}
\end{equation}

Further, for $\mathbb{n}$ large enough $\hat{\phi}_{t}$ is a Gaussian with mean $\phi_{t}$ and variance $\frac{1-\phi_{t}^2}{\mathbb{n}}$. Writing thus $\hat{\phi}_{t}=\phi_{t}+\eta_{t}$ where $\eta_{t}\sim \mathcal{N}\left(0,\frac{1-\phi_{t}^2}{\mathbb{n}}\right)$ and expanding Eq.~\ref{eq:2} at zeroth order, it is possible to write, for  $i=1, 2, \dots, N$:
\begin{equation}\label{eq:5}
\begin{aligned}
\lambda_{i,t+1} &=\left(\frac{\omega_i}{\lambda_{i,t}^2}+\frac{(1-\omega_i)\alpha^2\Sigma_\epsilon}{\left(1-\phi_{t}\right)^2}\right)^{-\frac{1}{2}}\doteq T_i(\{\lambda_{i,t}\}_{i=1,2,\dots,N})
\end{aligned}
\end{equation}
which is the evolution law already introduced in the previous section.
Note that, should one have considered the first order term in the expansion above, one would have obtained a heteroschedastic noise term superimposed on the deterministic skeleton given by the maps $T$. This kind of system has been studied in \cite{lillo2023analysis} in case $N=1$. Here, dealing with the more general case in which $N\geq1$, we aim to first understand the behavior of the deterministic skeleton before considering the random system that is based on it.

In the following, we will be interested in studying how the deterministic dynamical system defined by the $T$s behaves for different choices of the sets $\{\omega_i\}$ and $\{\pi_i\}$. 
One last remark must be done before proceeding. From the definition of leverage it is required that $\lambda_t\geq1\ \forall t\in \mathbb{N}$ and from the stationarity of the AR(1) process for $r_s$ it is required that $|\phi_t|\leq1$, i.e. $\sum_{i=1}^2\pi_i\lambda_{i,t}\leq 1+\gamma, \forall t\in\mathbb{N} $. When performing numerical simulations, initial conditions have been chosen at random in the box $[1,1+\gamma]^2$ and the simulations leading to violations of the constraints at any greater $t$ have been discarded.
In the following, let $\mathcal{V}\subseteq\mathbb{R}^2$ the largest set such that (for all the values of the $\omega$s and $\pi$s considered and for the given choices of the other parameters) by choosing initial conditions in this set the above constraints are satisfied. At least for $|\omega_1-\omega_2|$ sufficiently small, it can be shown that this set is nonempty and contains an open rectangle (see the supplementary material).


\subsection{Homogeneous case: sychronization}\label{Sec:EqMemories}
When the two banks use the same strategy in forecasting the risk, we obtain the following result.
\begin{theorem}{(Synchronization in the Homogenous Case)}\
If $\omega_1=\omega_2$, $\forall (\lambda_{1,0},\lambda_{2,0})\in\mathcal{V}$
\begin{equation}
\lim_{n\rightarrow\infty}|T_1^n(\lambda_{1,0},\lambda_{2,0})-T_2^n(\lambda_{1,0},\lambda_{2,0})|=0.\end{equation}
\end{theorem}
I.e. the leverages of different banks approach asymptotically the same orbit. This can be proven (see the supplementary material) by observing that the quantity $\frac{|\lambda_{1,t}-\lambda_{2,t}|}{\lambda_{1,t}+\lambda_{2,t}}$ is strictly decreasing along the orbits.
Therefore in the homogenous case it doesn't matter (as far as the behavior of the system is concerned) if a bank is larger than the other: the only relevant parameter to describe the evolution of the system is the common memory. This allows us to reduce this system to the one dimensional one already studied.
Lastly, this result can be easily generalized to the case in which $N>2$, obtaining that even if all the memories are the same then there will be asymptotic synchronization among the $N$ orbits regardless of the weights.
\subsection{Heterogeneous case}
Here the case in which $\omega_1\neq\omega_2$ is analyzed.

As mentioned, if $\pi_1\neq 0,1$ the behavior of the system depends nontrivially on $\pi_1,\omega_1$, and $\omega_2$ with the overall behaviour described by an interpolation of the behaviours associated to the isolated banks with memory parameters $\omega_1$ and $\omega_2$. This case won't be explored further here, as a complete rigorous understanding seems out of reach. Still, it is worth to emphasise that for some choices of the parameters a Hénon-like attractor seems to appear in the $\lambda_1$ vs $\lambda_2$ plane.

The box counting dimension for the attractor shown above has been calculated to be equal to $1.203\pm0.006$, thus confirming its apparently fractal nature (see the supplementary material for details).
Finally, the nature of the attractor led us to consider applying the main theorem in \cite{wang2008toward} about the existence of an SRB measure supported on it.
However, it seems that not all of the hypotheses can be met in our case: The main issue is that for $|\omega_1-\omega_2|$ going to zero, although the attractor becomes essentially one dimensional, there is no collapse of dimensionality for $\omega_1=\omega_2$ and initial conditions approach the diagonal only asymptotically; whereas,  for the H\'enon map, when the parameter usually denoted with ``$b$'' is zero there is a collapse of the dimensionality and all initial conditions are mapped to a 1D segment.


\subsubsection{'Big vs small' bank: the skew-product case}\label{Sec:Skew}
Let's focus here on the case with $\pi=0$. This case may be looked at as a forcing-forced system, in which the larger bank ``drives'' the smaller one.
This setting is indeed a subcase of the general ``heterogeneous case'' previously introduced, but the fact that several interesting results may be, even rigorously, obtained makes it an interesting case to study.

\begin{figure}
    \centering
    \includegraphics[width=\columnwidth]{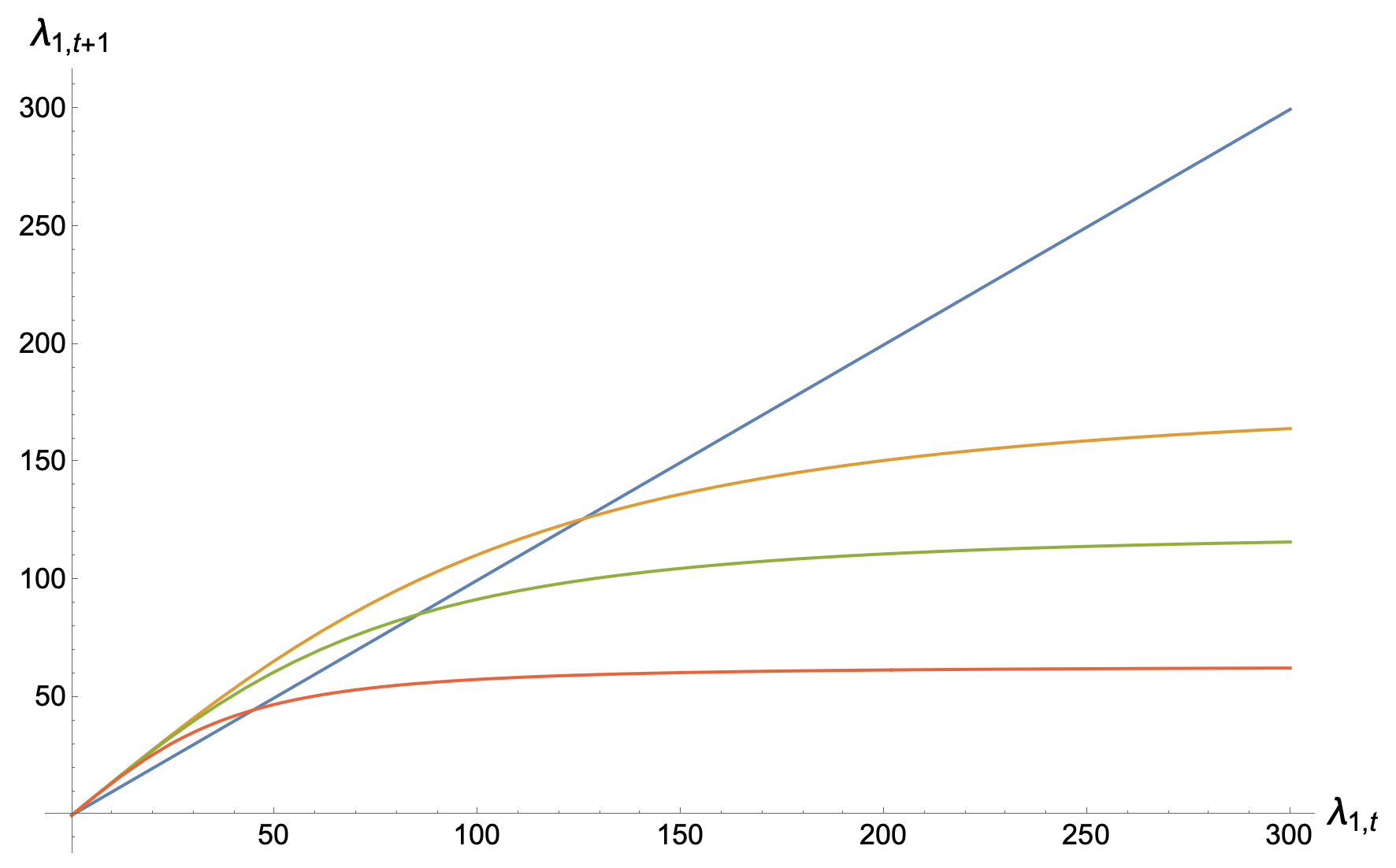}
    \caption{Some maps from the family of maps $f_y$ with $y=70$ (Orange), $y=80$ (Green), $y=90$ (Red). Here $\alpha=1.64,\Sigma_\epsilon=0.0015^2,\gamma=100,\omega_1=0.5$. The blue line is the bisector of the first quadrant.}
    \label{fig:mappefs}
\end{figure}

The dynamics of the leverage of the second bank (the forcing one) is described by the map $T$. The dynamics of the leverage of bank 1 depends  on the leverage of bank $2$. 
One could make this dependence explicit by introducing the family of maps $f_y$ so that \begin{equation}
    \lambda_{1,t+1}=f_{\lambda_{2,t}}(\lambda_{1,t})
\end{equation} with 
\begin{equation}
    f_{y}(x)=\left(\frac{\omega_1}{x^2}+\frac{(1-\omega_1)\gamma^2\alpha^2\Sigma_\epsilon}{(1+\gamma-y)^2}\right)^{-\frac{1}{2}}
\end{equation}
the maps $f_y$ are concave, monotonically increasing with a horizontal asymptote for large values of the argument (see Fig. \ref{fig:mappefs})

This system is in fact  a  skew product system
\begin{align}
\lambda_{1,t+1}&=f_{\lambda_{2,t}(\lambda_{1,t})}\\
\lambda_{2,t+1}&=T(\lambda_{2,t}).
\end{align}  To present some analytical results, it is useful to introduce the invertible extension of $T$, $\hat{T}:\hat{I}\rightarrow\hat{I}$ where \begin{equation}\hat{I}\doteq\{(\lambda_{2,i})_{i\in\mathbb{Z}}:T(\lambda_{2,i})=\lambda_{2,i+1}\forall i\in \mathbb{Z}\}\end{equation} and \begin{equation}\hat{T}((\lambda_{2,i})_{i\in\mathbb{Z}})=(\lambda_{2,i+1})_{i\in\mathbb{Z}}.\end{equation}
Also, for $\pmb{\lambda}\in\hat{I}$, let's write $f^n_{\pmb{\lambda}}=f_{\lambda_{2,n-1}}\circ f_{\lambda_{2,n-2}}\circ\dots\circ f_{\lambda_{2,1}}\circ f_{\lambda_{2,0}}$. (Let us stress once more that the evolution of the forcing bank does not depend on the forced bank and has already been discussed when dealing with the case $N=1$.)
Let's now give some analytical results (refer to supplementary material for further details and proofs).
To begin with, no matter the initial conditions, the leverages of the forced bank will always behave asymptotically in the same manner (determined by the initial condition of the large bank, as specified further in Theorem \ref{thm:fixed}). This also means that if one has two or more banks with weights $\pi_i=0$, their orbits will synchronize under the common forcing of the bigger bank.
\begin{theorem}{(Synchronization on the Fiber)}
For any orbit $\pmb{\lambda}\in \hat I$ of the forcing bank and for any initial conditions of the forced bank $\lambda_{1,0},\,\lambda_{1,0}^\prime\in[1,\infty)$ 
\begin{equation}
\lim_{n\rightarrow \infty}|\lambda_{1,n}-\lambda_{1,n}^\prime|=0.
\end{equation}
\end{theorem}

The proof of the theorem above is of topological nature and exploits the ``shape" of the functions $f_y$. In addition, one can also show that the Lyapunov Exponent of the forced bank is negative:
\begin{theorem}{(Negative Lyapunov Exponent on the Fiber)}
For any $\pmb{\lambda}\in \hat{I}$  and  $\lambda_{1,0}\in[1,\infty)$ 
\begin{equation}
\Lambda_1(\lambda_{1,0},\pmb{\lambda})\doteq\lim_{n\rightarrow\infty}\frac{1}{n}\log{|(f^n_{\pmb{\lambda}})^\prime(\lambda_{1,0})|}<0.
\end{equation}
\end{theorem}

From these results it looks like the dynamics of $\lambda_1$ is determined by something from the ``outside''. Not surprisingly, it is the large bank that determines the trajectory on which the leverage of the small bank will synchronize.
Let's start by looking at what happens when $\lambda_2$ reaches a periodic attractor

\begin{theorem}{(Periodic Forcing)}\\
    If the dynamics of the forcing bank is periodic, then the dynamics of the forced one will be periodic too (of the same period).
\end{theorem}
If instead the forcing is chaotic (meaning that it is topologically transitive on a finite union of intervals and that it has a unique a.c.i.m. with respect to Lebsegue), then the forced one shows some irregularity too:
\begin{theorem}{(Chaotic Forcing)}\\
If the dynamics of the forcing bank is chaotic, then the dynamics of the forced one is topologically transitive on an open subset of $[1,\infty)$. 
\end{theorem}
We emphasize that the proof of this theorem relies on the fact that the Schwarzian derivative of the map $T$ is negative, and that this is the unique assumption we have verified only numerically.
These results show how the choices of even a single parameter (e.g. the memory $\omega_2$) made by a large bank may have an impact on the stability of the whole system. This may be of interest to policymakers, who may for example constrain the values of the memory of the large bank in a ``safe'' region.

Next, it is possible to provide an explicit characterization of the behavior of the leverage of the small bank as a function of the past orbit of the forcing leverage. First, put, for $z\in I $:
\begin{equation}
a_{z}\doteq\frac{1+\gamma-z}{\gamma\alpha\sqrt{(1-\omega_1)\Sigma_\epsilon}}
\end{equation}

and define the function $x(\pmb{\lambda}): \hat{I}\rightarrow\mathbb{R}$
\begin{equation}
    x(\pmb{\lambda})\doteq\limsup_{n\rightarrow\infty}f^n_{\hat{T}^{-n}(\pmb{\lambda})}(a_{\lambda_2,-n})
\end{equation}
Where $\pmb{\lambda}=(\lambda_{2,i})_{i\in\mathbb{Z}}\in\hat{I}$.
Then one has:
\begin{theorem}{(Random Fixed Point)}\label{thm:fixed}\\ For every $\pmb{\lambda}\in\hat{I}$, it holds $f_{\lambda_{2,0}}(x(\pmb{\lambda}))=x(\hat{T}(\pmb{\lambda}))$ and
    \begin{equation}
        x(\pmb{\lambda})=\frac{1}{\sqrt{\sum_{i=0}^{\infty}A(\lambda_{2,-1-i})\omega_1^i}}
    \end{equation} with $A(\lambda)\doteq\frac{\gamma^2\alpha^2\Sigma_\epsilon}{(1+\gamma-\lambda)^2}$
\end{theorem}

Thus, for example, the fixed point depends 
continuously on the parameters and the past orbit of the forcing bank (e.g. with respect to the sup norm in the space $\hat{I}$).
In principle, this expression may for example allow the large bank to make the leverage of the small one behave in a specific way, provided it knows $\omega_1$. This is another evidence of the control the large bank has on the system.
Lastly, it is possible to make some estimates on the frequency of visit of subsets of the domain by the orbit of the forced bank by exploiting the ergodicity on the base (see the supplementary material). 

{\bf Supplementary Material}:
Proofs of theorems, as well as additional results and numerical simulations, are provided in the supplementary material file.

{\bf Acknowledgements}: The authors thank Fabrizio Lillo, Fulvio Corsi, and Sebastian van Strien for useful conversations. MT  acknowledges the support of EPSRC-FAPESP Grant No. 2023/13706 and EPSRC grant number UKRI1021. MT is very grateful to the Scuola Normale Superiore and the Centro di Ricerca Matematica Ennio De Giorgi for their hospitality during several visits in which this work was carried out.


\bibliography{references}

\end{document}